\documentclass{article}
\usepackage{graphicx}      

\usepackage[a4paper, total={6in, 8in}]{geometry}
\usepackage[numbers, sort]{natbib}
\usepackage[all]{xy}
\usepackage{caption}
\usepackage{mathrsfs}
\usepackage{placeins}
\usepackage{times} 
\usepackage{amsmath} 
\usepackage{amssymb}  

\usepackage{color}
\usepackage{graphics} 
\usepackage{epsfig} 
\usepackage{subfigure}
\usepackage{ulem}
\usepackage{tikz,pgfplots}
\usetikzlibrary{patterns}
\usepackage{algorithmic}
\usepackage{algorithm}

\newcommand{\N}{\mathbb{N}}

\newcommand{\R}{\mathbb{R}}
\newcommand{\grad}{\text{grad}}

\usepackage{bm}
\usepackage[colorlinks=true]{hyperref}
\hypersetup{urlcolor=blue, citecolor=red}
\newtheorem{theorem}{Theorem}[section]
\newtheorem{corollary}{Corollary}
\newtheorem{lemma}[theorem]{Lemma}
\newtheorem{proposition}{Proposition}

\newtheorem{remark}{Remark}

\usepackage{amssymb}

\newcommand{\proa}{A^*G \mbox{$\;$}_{\tau^*} \kern-3pt\times_\alpha
G \mbox{$\;$}_\beta \kern-3pt\times_{\tau^*} A^*G}

\tikzstyle{vertex}=[circle,fill=black!20,minimum size=15pt,inner sep=0pt]
\tikzstyle{selected vertex} = [vertex, fill=red!24]
\tikzstyle{edge} = [draw,thick,-]
\tikzstyle{dedge} = [draw,thick,<->]
\tikzstyle{shadowdedge} = [draw, dotted,->]
\tikzstyle{weight} = [font=\small]
\tikzstyle{selected edge} = [draw,line width=3pt,-,red!50]
\tikzstyle{ignored edge} = [draw,line width=3pt,-,black!20]

\title{Variational Obstacle Avoidance with Applications to Interpolation Problems in Hybrid Systems}
\author{Jacob R. Goodman and Leonardo J. Colombo}
\date{}

\begin{document}

\maketitle



\begin{abstract}
We study variational obstacle avoidance problems on complete Riemannian manifolds and apply the results to the construction of piecewise smooth curves interpolating a set of knot points in systems with impulse effects. We derive the dynamical equations for extrema in the variational problem, and show the existence of minimizers by using lower-continuity arguments for weak convergence on an infinite-dimensional Hilbert manifold. We then provide conditions under which it is possible to ensure that the extrema will safely avoid a given obstacle within some desired tolerance.
\end{abstract}


\section{Introduction}

Many problems in engineering, physics, biology, and related disciplines can be formulated as variational problems. A typical problem in this context is path planning. Sometimes the solution we seek has to satisfy some constraints or avoid static or moving obstacles in the space of configurations of a given system. It is also often the case that the desired paths must connect some set of knot points—interpolating positions with given velocities (and sometimes higher order derivatives too) \cite{CLACC, CroSil:95}. For such problems, the use of variationally defined curves has a rich history due to the regularity and optimal nature of the solutions. In particular, the so-called \textit{Riemannian splines} \cite{noakes} are a particularly ubiquitous choice in interpolant, which themselves are composed of Riemannian polynomials—satisfying boundary conditions in positions, velocities, and potentially higher-order derivatives—that are glued together. In Euclidean spaces, Riemannian splines are just cubic splines. That is, the minimizers of the total squared acceleration.

Riemannian polynomials are smooth and optimal in the sense that they minimize the average square magnitude of some higher-order derivative along the curve (a quantity which is often related to energy consumption in applications). Moreover, Riemannian polynomials carry a rich geometry with them, which has been studied extensively in the literature (see \cite{Giambo, marg, noakes, elastica} for a detailed account of Riemannian cubics and \cite{RiemannianPoly, popei} for some results with higher-order Riemannian polynomials). It is often the case that—in addition to interpolating points—there are obstacles or regions in space which need to be avoided. In this case, a typical strategy is to augment the action functional with an artificial potential term that grows large near the obstacles and small away from them (in that sense, the trajectories which minimize the action are expected to avoid the obstacles) \cite{kod}. This was done for instance in \cite{BlCaCoCDC} and \cite{BlCaCoIJC}, where necessary conditions for extrema in obstacle avoidance problems on Riemannian manifolds were derived, in addition to applications to interpolation problems on manifolds and to energy-minimum problems on Lie groups and symmetric spaces endowed with a bi-invariant metric. Nevertheless, there has been little to no work in the literature regarding safety guarantees and the role of potential shaping in successfully completing the task, which is the main focus of this paper. In particular, we will investigate the role of the artificial potential in the obstacle avoidance task on complete and connected Riemannian manifolds, and in doing so obtain conditions under which avoidance is guaranteed within some tolerance.

The main contributions of this paper are as follows: (1) We prove the existence of global minimizers to the variational problem in the case that the potential is smooth and non-negative, which is a necessary prerequisite in providing safety guarantees (indeed, proving that minimizing trajectories avoid an obstacle is useful only if such minimizing trajectories exist). This is accomplished by using standard techniques in functional analysis such as the weak lower semi-continuity of a norm in a Hilbert space, as was done in \cite{Giambo} in the case of Riemannian cubic polynomials. (2) We derive general conditions for the artificial potential—in terms of some reference trajectory which avoids the obstacle—under which the corresponding minimizers avoid a point-obstacle within some tolerance. We then remove the dependence on the reference trajectory for a particular family of potentials, and show that point-obstacle avoidance can be achieved within any desired tolerance for some potential in the family (constrained by the boundary conditions, the geometry of the manifold, and the sensing radius corresponding to technological limitations in the detecting the obstacle). (3) We consider the case where the obstacles are totally bounded subsets of the underlying manifold, and show that obstacle avoidance can still be achieved by treating the obstacle as a collection of point-obstacles and applying the previous techniques. (4) We apply these techniques to obtain natural and regular interpolants in multiple-domain hybrid systems, where Riemannian cubics splines may fail due to unintended collision with guards along the trajectory.


\section{Preliminaries on Riemannian Geometry}\label{Sec: background}


Let $Q$ be an $n$-dimensional  \textit{Riemannian
manifold} endowed with a non-degenerate symmetric covariant 2-tensor field $g$ called the \textit{Riemannian metric}. That is, to each point $q\in Q$ we assign an inner product $g_q:T_qQ\times T_qQ\to\mathbb{R}$, where $T_qQ$ is the \textit{tangent space} of $Q$ at $q$. The length of a tangent vector is determined by its norm,
$||v_q||=g(v_q,v_q)^{1/2}$ with $v_q\in T_qQ$. A \textit{Riemannian connection} $\nabla$ on $Q$ is a map that assigns to any two smooth vector fields $X$ and $Y$ on $Q$ a new vector field, $\nabla_{X}Y$. For the properties of $\nabla$, see \cite{Boothby}.  The operator
$\nabla_{X}$, which assigns to every vector field $Y$ the vector
field $\nabla_{X}Y$, is called the \textit{covariant derivative of
$Y$ with respect to $X$}.


Consider a vector field $W$  along a curve $q$ on $Q$. The $k$th-order covariant derivative  of $W$ along $q$ is denoted by $\displaystyle{\frac{D^{k}W}{dt^{k}}}$, $k\geq 1$. We also denote by $\displaystyle{\frac{D^{k+1}q}{dt^{k+1}}}$ the $k$th-order covariant derivative of the velocity vector field of $q$  along $q$, $k\geq 1$.

A vector field $X$ along a piecewise smooth curve $q$ in $Q$ is said to be \textit{parallel along $q$} if $\displaystyle{\frac{DX}{dt}\equiv 0}$.
Given vector fields $X$, $Y$
and $Z$ on $Q$, the vector field $R(X,Y)Z$ given by $
R(X,Y)Z=\nabla_{X}\nabla_{Y}Z-\nabla_{Y}\nabla_{X}Z-\nabla_{[X,Y]}Z$ is called the \textit{curvature endomorphism} on $Q$. $R$ is trilinear in $X$, $Y$ and $Z$.


If we assume that $Q$ is \textit{complete}, then any two points $x$ and $y$ in $Q$ can be connected by a minimal length geodesic $\gamma_{x, y}$, and the Riemannian distance $d:Q\times Q\to\mathbb{R}$ between two points in $Q$ can be defined by
$\displaystyle{d(x,y)=\int_{0}^{1}\Big{\|}\frac{d \gamma_{x,y}}{d s}(s)\Big{\|}\, ds}$. The idea of a geodesic is useful because it provides a map from $T_qQ\hbox{ to }Q$ in the following way: $v\in T_{q}Q \mapsto \gamma(1),\, \gamma(0)=q,\, \dot{\gamma}(0)=v$, where $\gamma$ is a geodesic. This map is called the \textit{Riemannian exponential map} and is denoted by $\mathrm{exp}_q:T_qQ\to Q$. In particular, $\mathrm{exp}_q$ is a diffeomorphism from some star-shaped neighborhood of $0 \in T_q Q$ to a \textit{geodesically convex} open neighborhood $\mathcal{B}$ of $q \in Q$. That is, any two points in $\mathcal{B}$ can be connected by a unique minimizing geodesic. Moreover, if $y \in \mathcal{B}$, we can express the Riemannian distance locally by means of the Riemannian exponential as $d(q,y)=\|\mbox{exp}_q^{-1}y\|.$

\vspace{.2cm}

The \textit{Lebesgue space} $L^p([0,1];\R^n)$, $p\in(1,+\infty)$ is the space of $\R^n$-valued functions on $[0,1]$ such that each of their components is $p$-integrable, that is, whose integral of the absolute value raised to the power of $p$ is finite. 
A sequence $(f_n)$ of functions in $L^p([0,1];\R^n)$ is said to be \textit{weakly convergent} to $f$ if for every $g\in L^r([0,1];\R^n)$, with  $\frac{1}{p}+\frac{1}{r}=1$, and every component $i$,  $\displaystyle{\lim_{n\to\infty}\int_{[0,1]}f_n^i g^i=\int_{[0,1]}f^i g^i}$. A function $g\colon[0,1]\to \R^n$ is said to be the \textit{weak derivative} of $f\colon[0,1]\to \R^n$ if for every component $i$ of $f$ and $g$, and for every compactly supported $\mathcal{C}^\infty$ real-valued function $\varphi$ on $[0,1]$, $\displaystyle{\int_{[0,1]}f^i\varphi'=-\int_{[0,1]}g^i\varphi}$.  The \textit{Sobolev space} $W^{k,p}([0,1];\R^n)$ is the space of functions $u\in L^p([0,1];\R^n)$ such that for every $\alpha\leq k$, the $\alpha^{th}$ weak derivative $\frac{d^\alpha u}{dt^{\alpha}}$ of $u$ exists and $\frac{d^\alpha u}{dt^{\alpha}}\in L^p([0,1];\R^n)$. In particular, $H^k([0,1];\R^n)$ denotes the Sobolev space $W^{k,2}([0,1];\R^n)$, and its norm may be expressed as $\displaystyle{\left|\left|f \right| \right| = \left(\int_{[0,1]} \sum_{p=0}^k \left| \left|\frac{d^k}{dt^k}f(t)\right|\right|_{\R^n}^2 dt \right)^{1/2}}$ for all $f \in H^k([0,1];\R^n)$, where $|| \cdot ||_{\R^n}$ denotes the Euclidean norm on $\R^n$. $(f_n)\subset W^{k,p}([0,1];\R^n)$ is said to be \textit{weakly convergent} to $f$ in $W^{k,p}([0,1];\R^n)$ if for every $\alpha\leq k$, $\displaystyle{\frac{d^\alpha f_n}{dt^{\alpha}}\rightharpoonup \frac{d^\alpha f}{dt^{\alpha}}}$ weakly in $L^p([0,1];\R^n)$.


We denote by $H^2([0,1];Q)$ the set of all curves $q\colon[0,1]\to Q$ such that for every chart $(\mathcal{U},\varphi)$ of $Q$ and every closed subinterval $I\subset[0,1]$ such that $q(I)\subset\mathcal{U}$, the restriction of the composition $\varphi\circ q|_I$ is in $H^2([0,1];\R^m)$. Note that $H^2([0,1]; Q)$ is an infinite-dimensional Hilbert Manifold modeled on $H^2([0,1]; \R^m)$, and given $\xi = (q_0,v_0), \ \eta = (q_T, v_T) \in TQ$, the space $\Omega_{\xi, \eta}^T$ (denoted simply by $\Omega$ unless otherwise necessary) defined as the space of all curves $\gamma \in H^2([0,1]; Q)$ satisfying $\gamma(0) = q_0, \ \gamma(T) = q_T, \ \dot{\gamma}(0) = v_0, \ \dot{\gamma}(T) = v_T$ is a closed submanifold of $H^2([0,1]; Q)$ (see \cite{Palais}). The tangent space $T_x \Omega$ consists of vector fields along $x$ of class $H^2$ which vanish at the endpoints together with their first covariant derivatives. We consider the Hilbert structure on $T_x \Omega$ induced by the inner product $\displaystyle{\left< V, W \right> := \int_0^T g\left(\frac{D^2}{dt^2}V, \frac{D^2}{dt^2} W\right)dt}$. This inner product induces (fiberwise) a Riemannian metric on $\Omega$, which itself induces a metric in the usual way. It is known that the completeness of $\Omega$ follows from the completeness of $Q$ (see \cite{hilbert}).

\section{Variational Obstacle Avoidance Problem}\label{Sec: Necessary conditions}


Consider a complete connected Riemannian manifold $Q$, and the space $\Omega$ as defined in Section \ref{Sec: background}. We define $J: \Omega \to \R$ as
\begin{equation}\label{J}
J(q)=\frac{1}{2} \int\limits_0^T \Big{(}\Big{|}\Big{|}\frac{D^2q}{dt^2}(t)\Big{|}\Big{|}^2 + V(q(t))\Big{)}dt.
\end{equation}

\textbf{Variational obstacle avoidance problem:} Find a curve $q\in \Omega$ minimizing the functional $J$, where $V:Q \to\mathbb{R}$ is a smooth and non-negative function called the \textit{artificial potential}.

\vspace{.2cm}

In order to minimize the functional $J$ among the set $\Omega$, we want to find curves $q\in\Omega$ such that $J(q)\leq J(\tilde{q})$ for all
admissible curves $\tilde{q}$ in an $H^2$-neighborhood of $q$. The next result from \cite{BlCaCoCDC} characterizes necessary conditions for optimality in the variational obstacle avoidance problem.


\begin{proposition}\label{th1}\cite{BlCaCoCDC} A curve $q \in \Omega$ is a critical point of the functional $J$ if and only if it is smooth on $[0,T]$ and satisfies:
\begin{equation}\label{eqq1}
    \frac{D^4q}{dt^4}+R\Big{(}\frac{D^2q}{dt^2},\frac{dq}{dt}\Big{)}\frac{dq}{dt}=- \hbox{\grad} \, V(q(t)).
\end{equation}
\end{proposition}

\begin{remark}
In \eqref{J}, the artificial potential $V$ is introduced for the purpose of obstacle avoidance. Repulsive potentials (as we will study in Section \ref{Sec: Obstacle Avoidance}) are particularly well-suited for this goal. However, this is not integral to the variational problem considered (indeed, smooth and non-negative potentials encompasses a very broad range of functions), so that the above formalism is adaptable to variational problems with collective behavior performances other than obstacle avoidance—such as collision avoidance of multi-agent systems \cite{sh, CoGo20, CollAvoid}. 
\end{remark}

\section{Existence of global minimizers}\label{Sec: existence}

Next, we will prove the existence of global minimizers of $J$ in $\Omega$ by employing some classical techniques from functional analysis, such as lower semi-continuity arguments for weak convergence in $H^2$. Before beginning the proof we will introduce a lemma that simplifies the analysis considerably.

\begin{lemma}\label{Convergence}
Let $Q$ be an $m$-dimensional complete Riemannian manifold, and suppose that $\{q_n\} \subset \Omega$ is a sequence such that $\displaystyle{\sup_{n \in \N} J(q_n) < +\infty}$. Then, $\{q_n\}$ and $\{\dot{q}_n\}$ are uniformly bounded, and there exists a subsequence of $\{q_n\}$ which converges weakly to some $q \in \Omega$ with respect to the norm on $H^2$.
\end{lemma}

\textit{Proof:} This follows from Lemma 4.2 of \cite{CollAvoid} in the case of $s = 2$, with one agent being static.

\vspace{.2cm}

\begin{theorem}\label{Existence}
The functional $J$ attains its minimum in $\Omega$.
\end{theorem}
\textit{Proof:} Suppose that $\{q_n\} \subset \Omega$ is a minimizing sequence.
That is, $\lim_{n \in \N} J(q_n) = \inf_{q \in \Omega} J(q) \ge 0$. Note that, such a sequence satisfies the assumptions of Lemma \eqref{Convergence}, so that there exists a subsequence of $q_n$ (which we also denote by $\{q_n\}$ for convenience) that converges weakly to some $q \in \Omega$ with respect to the norm on $H^2$. It then suffices to show that $\displaystyle{J(q) \le \liminf_{n \to \infty} J(q_n)}$. Since $g(\dot{q}_n, \dot{q}_n)$ is uniformly bounded, there exists a finite collection of charts $(U_i, \varphi_i)$ on $Q$ and $I_i$ an accompanying finite partition of $[0,1]$ such that, for sufficiently large $n$, there exists a compact subset $K_i \subset U_i$ containing $q_n(I_i)$. In local coordinates, we may consider $q_n$ to be a curve on $\R^m$ (however, we will abuse this notation by continuing to call it $q_n$ both on the chart $U_i$ and its image in $\R^m$).

Observe first that $V(q_n)$ converges to $V(q)$ uniformly on $[0, T]$ since the interval is compact, $V$ is continuous, and $q_n \to q$ uniformly. Therefore, showing that $J(q) \le$ $\displaystyle{\liminf_{n \to \infty} J(q_n)}$ is then equivalent to showing that for all intervals $\mathcal{I}_i$, we have
\begin{small}
\begin{align*}
\int_{\mathcal{I}_i} g \left(\frac{D}{dt} \dot{q}(t), \frac{D}{dt} \dot{q}(t)\right)dt &\le \liminf_{n \to \infty} \int_{\mathcal{I}_i} g \left(\frac{D}{dt} \dot{q}_n(t), \frac{D}{dt} \dot{q}_n(t)\right)dt.
\end{align*}
\end{small}
Note that $\frac{d}{dt}\dot{q}_n = \frac{D}{dt} \dot{q}_n + \Gamma(q_n; \dot{q}_n, \dot{q}_n)$ where $\Gamma: \R^{3m} \to \R^m$ is continuous in the first argument and bilinear in the last two—and it is determined by the induced Christoffel symbols. It follows that $\Gamma(q_n; \dot{q}_n, \dot{q}_n) \to \Gamma(q; \dot{q}, \dot{q})$ uniformly on $I_i$, so that the above inequality is equivalent to
\begin{align}\label{riem_ineq}
     \int_{\mathcal{I}_i} g \left( \ddot{q}(t), \ddot{q}(t) \right)dt   &\le \liminf_{n \to \infty} \int_{\mathcal{I}_i} g \left(\ddot{q}_n(t), \ddot{q}_n(t)\right)dt.
\end{align}
Thus far, we have suppressed the dependence of the Riemannian metric on the point at which we are evaluating the tangent vectors. Note that this is not problematic by the uniform convergence of $q_n$ to $q$. That is, 
\begin{align*}
    \left|\int_{\mathcal{I}_i} g_{q(t)} \left(\ddot{q}_n(t), \ddot{q}_n(t)\right)dt - \int_{\mathcal{I}_i} g_{q_n(t)} \left(\ddot{q}_n(t), \ddot{q}_n(t)\right)dt \right| \to 0,
\end{align*} as $n \to \infty$, so we may assume that the metric is evaluated at $q(t)$ on both sides of inequality \eqref{riem_ineq}. We now consider the set 
\begin{small}
\begin{align*}
     L_g^2(I_i, \R^m) := \Big{\{} \gamma: I_i \to \R^m \ : \int_{\mathcal{I}_i} g_{q(t)} \left( \gamma(t), \gamma(t) \right)dt < +\infty \Big{\}},
\end{align*}
\end{small}
which can be endowed with the structure of a normed linear space, with the norm $\displaystyle{||\gamma|| = \int_{\mathcal{I}_i} g_{q(t)} \left( \gamma(t), \gamma(t) \right)dt}$.

Note that the Euclidean and Riemannian norms are (bi-Lipchitz) equivalent in the compact chart image of $TQ |_{K_i}$, which further implies that $L^2(I_i, \R^m)$ and $L^2_g(I_i, \R^m)$ are equivalent as normed linear spaces (indeed, $L^2_g$ can be thought of as a weighted $L^2$ space in local coordinates). Hence, they induce the same weak topology, and so the weak $L^2$-convergence of $\ddot{q}_n$ to $\ddot{q}$ further implies its weak $L^2_g$-convergence—from which \eqref{riem_ineq} follows immediately. $\hfill\square$

\section{The Obstacle Avoidance Task}\label{Sec: Obstacle Avoidance}
In this section, we explore the task of obstacle avoidance. In particular, we derive conditions on the artificial potential under which obstacle avoidance is guaranteed within some tolerance. Section \ref{sec: Point-obstacles} handles the case that the obstacle is a point on the manifold $Q$. We then extend the analysis in Section \ref{sec: large-obstacles} to the case where the obstacle is a totally bounded subset of $Q$.


\subsection{Point-Obstacles}\label{sec: Point-obstacles}

We first consider the case of a point-obstacle. We fix some point $p \in Q$ as the obstacle. Denote by $B_s(p)$ the ball of radius $s$ centered at $p.$ Let $0 < r < r^\ast < R$, and define the \textit{Collision region} $C_p := B_r(p)$, the \textit{Risk region} $C^\ast_p := B_{r^\ast}(p)$, and the \textit{Safety region} $S_p := Q \setminus \overline{B_R}(p)$, where $\overline{B_R}(p)$ denotes the topological closure of $B_R(p)$. We say that $q \in \Omega$ avoids the obstacle $p$ with tolerance $r$ if $q(t) \notin C_p$ for all $t \in [0, T]$. We will construct our potential so that it is bounded above by some constant in the Safety regions and bounded below by some (larger) constant in the Risk region. More precisely, for some real numbers $0 \le V^- \le V^\ast$, we construct the artificial potential such that $V \ge V^\ast$ on $C^\ast_{p}$ and $V \le V^{-}$ on $S_{p}$.

We call obstacle avoidance with tolerance $r$ \textit{feasible} if there exists a curve $q(t) \in \Omega$ avoiding the obstacle with tolerance $r$.  We further define a \textit{reference trajectory} $q \in \Omega$ to be a curve $q(t) \in S$ for all $t \in [0, T]$ (note that this requires $R < \min\left\{d(p, q_0), \ d(p, q_T)\right\}$). Clearly, the existence of a reference trajectory is equivalent to saying that obstacle avoidance is feasible with the tolerance $R$ of the Safety region.

\begin{proposition}\label{Collision_Avoidance_Prop}
Assume collision avoidance is feasible with the tolerance $R$, and fix a reference trajectory $q \in \Omega$. Define the non-negative real numbers $\displaystyle{a := \sup_{t \in [0, T]} \left|\left|\frac{D}{dt}\dot{q}(t)\right|\right|}$, $c := (a^2 + V^-)T$, and $v:=\sqrt{cT} + \sqrt{cT + ||v_0||^2}$. If $V^\ast > \frac{cv}{2(r^\ast - r)}$, then any minimizer $q^\ast$ of $J$ avoids $p$ with tolerance $r$.
\end{proposition}

\textit{Proof:} 
It is clear that the variational collision avoidance problem (as defined in \cite{sh, CoGo20, CollAvoid}) reduces to the variational obstacle avoidance problem in the case of two agents, with the second agent constrained as $q_2(t) \equiv p$ for all $t \in [0, T]$. From here, the proof follows immediately by Proposition 2 in \cite{CollAvoid}.

\begin{remark}\label{remark: sensing}
In applications, distances can be calculated by, for instance, attaching a sensor to the agent. It is often the case that—due to technological limitations—measurements on distance are only possible/reliable within some given ball of the sensor, the radius $h$ of which is called the \textit{sensing radius}. We may account for this by demanding that the potential vanish identically whenever the distance between the sensor and the obstacle exceeds $h$. This can be handled—while still preserving the regularity and positive-definiteness of the obstacle avoidance potential—by utilizing bump functions. Observe that this does not affect the analysis of Proposition \ref{Collision_Avoidance_Prop}. 
\end{remark}

Accounting for the sensing radius $h$ as discussed in Remark \ref{remark: sensing}, we now consider the smooth, non-negative family of potentials parameterized by $D,\tau \in \R^+, \ k \in \N$ defined by $$V_{D, \tau}^k(q) = \begin{cases} e\tau\exp\left(-\frac1{1 - (d(p, q)/D)^{2k}}\right) & d(p, q) < D \\
0 & \text{else}
\end{cases}$$

Observe that when $D \le h$, we satisfy the requirement that the potential vanish identically outside of the sensing radius.

\begin{corollary}\label{cor: point_obstacle_avoidance}
If obstacle avoidance is feasible for the tolerance $0 < R \le h$, then for all $r < R$, there exists $\tau^\ast \in \R^+, \ k^\ast \in \N$ such that for all $\tau > \tau^\ast$ and $k > k^\ast$, every minimizer of $J$ with the potential $V = V_{R, \tau}^{k}$ avoids the obstacle with tolerance $r$.
\end{corollary}
\textit{Proof:} Since obstacle avoidance is feasible with tolerance $R$, there exists a reference trajectory $q \in \Omega$ that remains in the safety region with tolerance $R$ for all $t \in [0, T]$. Fix $r < R$ and choose $r^\ast$ such that $r < r^\ast < R$. It follows by definition that $V_{R, \tau}^k(q) = 0$ whenever $d(p,q) \ge R$—independent of $\tau, k$. Hence we have $V^- = 0$.

Moreover, $d(p, q)/R < r^\ast/R < 1$ whenever $d(p, q) < r^\ast$, so that $(d(p, q)/R)^{2k}$ can be made arbitrarily small by taking $k$ sufficiently large. Since $V_{R, \tau}^k$ is continuous, it follows that for all $\tau, \epsilon > 0$, there exists a $K > 0$ such that for all $k > K$, we have $V_{R, \tau} > \tau - \epsilon$ on $C_p^*$. In other words, $V^\ast$ can be made arbitrarily large by choosing the parameters appropriately. Since $\frac{cv}{2(r^\ast - r)}$ is finite and independent of the parameters $\tau, k, D$, it follows by Proposition \ref{Collision_Avoidance_Prop} that any minimizer of $J$ avoids the obstacle with tolerance $r$ by choosing $\tau$ (and hence $k$) sufficiently large. $\hfill\square$

\subsection{General Obstacles}\label{sec: large-obstacles}

We now focus our attention to avoiding larger obstacles, which in the most general case are simply modeled as subsets of $Q$ with no additional structure. Similar to the case of point-obstacles, we will say that the obstacle avoidance task (with obstacle $P \subset Q$) is feasible with tolerance $r$ if there exists some $q \in \Omega$ such that $d(q(t), P) > r$ for all $t \in [0, T]$, where for $m \in Q$, the distance function is defined as $d(m, P) := \inf_{p \in P} d(m, p)$. For points $m, p \in Q$, it is easy to see that the triangle inequality holds as $d(m, P) \le d(m, p) + d(p, P)$. From this it is clear that the analysis carried out in Proposition \ref{Collision_Avoidance_Prop} and  Corollary \ref{cor: point_obstacle_avoidance} follows identically, with all occurrences of the distance from $p \in Q$ replaced by the distance from $P \subset Q$. 

\begin{corollary}\label{cor: Obstacle_avoidance}
If obstacle avoidance is feasible for the tolerance $0 < R \le h$, then for all $r < R$, there exists $\tau^\ast \in \R^+, k^\ast \in \N$ such that for all $\tau > \tau^\ast$ and $k > k^\ast$, every minimizer of $J$ with $$V_{R, \tau}^k(q) = \begin{cases} e\tau\exp\left(-\frac1{1 - (d(P, q)/R)^{2k}}\right) & d(P, q) < R \\
0 & \text{else},
\end{cases}$$ avoids the $P$ with tolerance $r$.
\end{corollary}

This corollary poses an issue for practical purposes, as calculating $d(P, q)$ is computationally intensive. Along a minimizer, we would in general need to re-calculate the distance for each point in time. Moreover, even in the case where the boundary points $q_0, q_T \in Q$ for $\Omega$ are sufficiently close to represent the distance between points via the Riemannian exponential map as $d(m, p) = \left|\left| \exp_m^{-1}(p)\right|\right|$, the gradient $\grad_1 d(m, P)$ may not have a closed form (and so neither will $\grad V(m)$). We therefore seek an alternative characterization which relies only on the distance function and family of potentials used in the previous section. In particular, we will suppose that the obstacle $P$ is totally bounded, so that we may cover it in a finite collection of open balls—the centers of which will be treated as point-obstacles. The potential will then be taken as a sum of the point-obstacle avoidance potentials that result in the minimizers of $J$ avoiding these balls. Before stating the result, we introduce a lemma describing some basic facts of totally bounded sets that will be useful in the subsequent proof.

\begin{lemma}\label{lemma: totally_bounded}
Suppose $P \subset Q$ is totally bounded. Then,
\begin{enumerate}
    \item $\overline{P}$ is totally bounded.
    \item Any subset $A \subset P$ is totally bounded.
    \item $B_r(P)$ is totally bounded for any $r > 0$.
\end{enumerate}
\end{lemma}

\begin{proposition}\label{prop: obst_avoid}
Suppose that obstacle avoidance is feasible with tolerance $0 < R \le h$. Then for all $r < \frac{R}2$, there exists an $r^\ast < R$, $\tau > 0$, $k \in \N$, and a finite collection of points $\{p_i\} \subset P$ such that any minimizer of $J$ corresponding to the potential
$V(q) = \sum_i V_i(q)$, where $V_i(q)$ is defined by
$$V_i(q) = \begin{cases} e \tau\exp\left(-\frac1{1 - (d(p_i, q)/R)^{2k}}\right) & d(p_i, q) < R \\
0 & \text{else},
\end{cases}$$
avoids $P$ with tolerance $r$. 
\end{proposition}
\textit{Proof:} Fix $r < \frac{R}{2}$. By lemma \ref{lemma: totally_bounded}, the set $\overline{B_r}(P)$ is totally bounded and hence $\partial B_r(P) \subset \overline{B_r}(P)$ is totally bounded, where $\partial B_r(P)$ denotes the boundary of $B_r(P)$. Define the real numbers $\delta := \frac{R - 2r}{2}$ and $r^\ast := \frac{2r + \delta + R}{2}$ so that $0 < r < 2r + \delta < r^\ast < R$. Now let $\{p_i^\ast\} \subset \partial B_r(P)$ be a finite collection of points such that $\partial B_r(P) \subset \cup_i B_{\delta}(p_i^\ast)$. For each $p_i^\ast$, choose a point $p_i \in P$ such that $d(p_i, p_i^\ast) = r$ (such a point exists by the definition of $\partial B_r(P)$). 

Next we show that the collection of balls $\{B_{r^\ast}(p_i)\}$ cover $B_r(P)$. To that end, let $m \in B_r(P)$, $m^\ast$ be a point on $\partial B_r(P)$ which minimizes the distance from $m$, and $p_j^\ast$ be an element of $\{p_i^\ast\}$ which minimizes the distance from $m^\ast$. Observe that $ d(m, p_j) \le d(m, p_j^\ast) + d(p_j^\ast, p_j) \le d(m, m^\ast) + d(m^\ast, p_j^\ast) + r \le 2r + \delta < r^\ast$, so that $m \in B_{r^\ast}(p_j)$. From this and the fact that $r^\ast < R$, we have $B_r(P) \subset \cup_i B_{r^\ast}(p_i) \subset B_R(P)$. From Corollary \ref{cor: point_obstacle_avoidance}, it is clear that we may make the potential $V_i$ arbitrarily large on $B_{r^\ast}(p_i)$ by choosing $\tau, k$ sufficiently large. Since every potential of this form is positive-definite, it follows that we may make $V = \sum_i V_i$ arbitrarily large on each $B_{r^\ast}(p_i)$, and hence on $B_r(P)$. Moreover, on $Q \setminus \overline{B_R}(P)$, we have $V \equiv 0$. The remainder of the proof follows from feasibility of the obstacle avoidance task and similar analysis to that found in the proofs of Proposition \ref{Collision_Avoidance_Prop} and Corollary \ref{cor: point_obstacle_avoidance}. 
$\hfill\square$

\begin{remark}
Proposition \ref{prop: obst_avoid} implies that we can avoid a finite collection of point-obstacles (when feasible) by summing over an appropriate family of potentials—each of which correspond to the avoidance of one of the point-obstacles. Similarly, we can avoid a finite collection of totally bounded obstacles (when feasible) by summing over appropriate potentials.
\end{remark}

\subsection{Simulation Results}
In this section, we conduct numerical simulations of the collision avoidance problem in the case of $Q = \R^3$ with the Euclidean metric. We choose our obstacle as a section of the unit sphere. In particular, we choose $P = \{(\sin(\phi)\sin(\theta), \sin(\phi)\cos(\theta), \cos(\phi)) \in \R^3 \ : \ 0 < \phi < \pi/4, \ 0 < \theta < \pi/2 \}$. Boundary conditions were chosen so that solutions to the boundary value problem with the potential $V \equiv 0$ collide with the obstacle. In particular, we choose:
\begin{align*}
    &q(0) = (0, 0, 0), &&\dot{q}(0) = (0.125, 0.125, 0.45) \\
    &q(1) = (0.2, 0.5, 1.8), &&\dot{q}(1) = (0.3, 0.25, 0.5)
\end{align*}
 
Motivated by Section \ref{sec: large-obstacles}, we choose three point obstacles $p_1, p_2, p_3$ of the form \\$(x, y, z) = (\sin(\phi)\sin(\theta), \sin(\phi)\cos(\theta), \cos(\phi))$ with $(\phi, \theta) = (\frac{\pi}{12}, \frac{\pi}{4}), \ (\frac{\pi}{5}, \frac{\pi}{9}), \ (\frac{\pi}{5}, \frac{\pi}{3})$ so that the balls of radius $R = 0.3$ centered at these points completely cover the obstacle $P$. We then construct the potential $V$ as in Proposition \ref{prop: obst_avoid}, namely so that we avoid the points $p_1, p_2, p_3$ with the tolerance of $R = 0.3$—and hence avoid $P$—with $k = 4$ and $\tau = 100/e$.

\begin{figure}[h!]
\begin{center}
 \includegraphics[width=6cm]{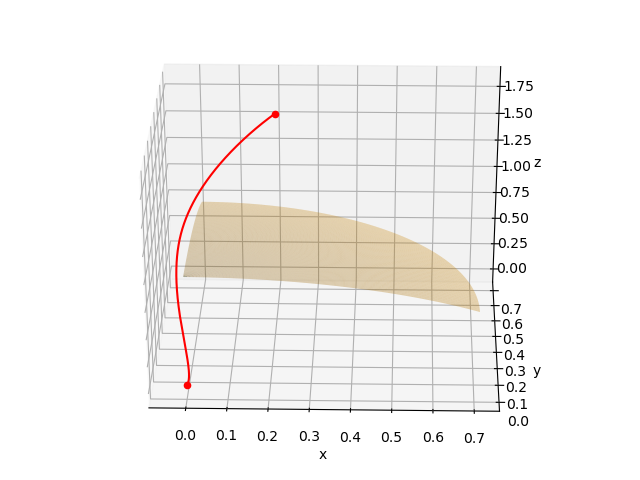}
 \includegraphics[width=6cm]{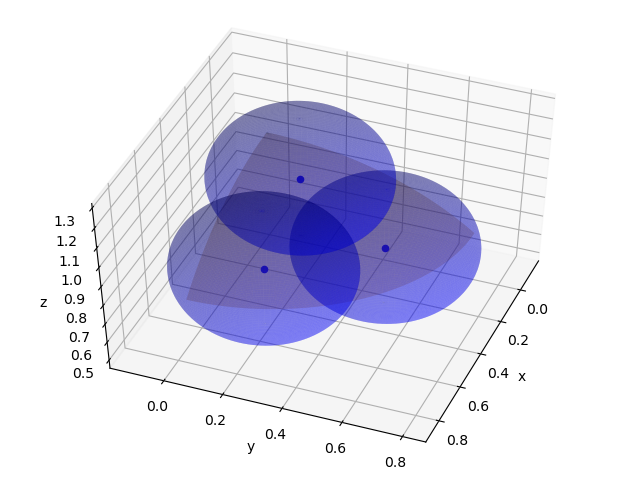}
 \caption{The plot on the left shows the numerical solutions to equations \eqref{eqq1} with the potential $V$ as defined above, from which it is clear that the agent avoids the obstacle $P$ while solving the boundary value problem. The plot on the right shows the covering of the obstacle $P$ by balls of radius $R = 0.3$ centered at the points $p_1, p_2,$ and $p_3$.}\label{fig: R3}
 \label{figuav}
 \end{center}
\end{figure}

In Figure \ref{fig: R3}, we plot the solutions to the necessary conditions \eqref{eqq1} with the potential $V$ as defined above. Numerical integration is done via the Euler method with a time step of $h = 0.001$. A shooting method based on the downhill simplex algorithm was used to find the initial accelerations and jerks that lead to solutions to the boundary value problem.

\section{Interpolation in Systems with Impulses}\label{sec: orbits}

\subsection{Multiple-Domain Hybrid Dynamical Systems}

A \textit{Dynamical system with impulse effects} (SIEs) is a class of dynamical system that exhibits both discrete and continuous behaviors \cite{teelSurvey}, \cite{CoGo19}, \cite{cha}. The transition from one to the other occurs when the continuous-time flow reaches a co-dimension one submanifold of the state space and re-initializes via a $C^1$ map on this submanifold.

\textit{Definition:} A $k$-domain SIEs is a tuple $\mathscr{H} = (\Gamma, Q, S, \Delta, X)$, where:
\begin{itemize}
\item[(i)] $\Gamma = (\mathcal{V}, \mathcal{E})$ is a connected, directed graph such that $\mathcal{V} = \{y_1,..., y_k\}$ is a set of $k$ vertices, and $\mathcal{E} \subset Q \times Q$ is the  set of edges. We further define the set of neighbors of a vertex $i \in \mathcal{V}$ by $\mathcal{N}_i = \{ j \in \mathcal{V}: e_{ij} \in \mathcal{E} \}$.
\item[(ii)] $Q = \{Q_y\}_{y \in \mathcal{V}}$ is a collection of smooth complete, connected Riemannian manifolds with Riemannian metric $g_y$.
\item[(iii)] $S = \{S_e\}_{e \in \mathcal{E}}$ is a collection of guards, where for $e = (y_i, y_j) \in \mathcal{E}$, $S_e$ is assumed to be a totally bounded embedded open section of $Q_{y_i}$.
\item[(iv)] $\Delta = \{\Delta_e\}_{e \in \mathcal{E}}$ is a collection of reset maps, which are $C^1$ mappings where, for $e = (y_i, y_j) \in \mathcal{E}$, \ $\Delta_e: S_e \to M_{y_j}$.
\item[(v)] $X = \{X_y\}_{y \in \mathcal{V}}$ is a collection of smooth vector fields $X_y: Q_y \to TQ_y$.
\end{itemize}
The underlying dynamical system with impulse effects is then defined by
$$\begin{cases} 
      \dot{q} = X_{i}(q) & \ \text{ if } \ q \in Q_{i} \setminus \cup_{j \in \mathcal{N}_i} S_{{ij}} \\
      q^+ = \Delta_{ij}(q^-) & \ \text{ if } \ q^- \in S_{ij} \ \text{for some} \ j \in \mathcal{N}_i,
\end{cases}$$
where it is understood that $X_i = X_{y_i}$, $Q_i = Q_{y_i}$, and $S_{ij} = S_{e_{ij}}$, $\Delta_{ij} = \Delta_{e_{ij}}$, and where $q^-$ and $q^+$ denote the left-limit and right-limit, respectively, of the trajectory as it intersects $S_i$ (and is correspondingly reset by $\Delta$). In general, $q^- \ne q^+$, so that there may be a point of discontinuity here. However, as in \cite{CoGo19}, we are given the choice in deciding whether the trajectory will be left-continuous or right-continuous at this point. That is, whether $q^- \in S_{ij}$ or $q^+ \in \Delta_{ij}(S_{ij})$ belong to our trajectory. In this paper, we will choose the former. Note that the results that follow in this work hold regardless of this choice \cite{cha}.

\begin{remark}\label{zeno} The above definition may in principle lead to a phenomenon known as \textit{Zeno behavior}, characterized by an infinite number of resets in finite time—which is particularly problematic in applications where numerical work is used, as computation time grows infinitely large at these Zeno points. There are two primary modes through which Zeno behavior can occur:

\begin{enumerate}
    \item A trajectory is continuously reset back to a guard. To exclude this type of behavior, we require that \\ $\Big{(} \cup_{k \in \mathcal{N}_j} S_{jk} \Big{)} \cap \overline{\Delta_{ij}}(S_{ij}) = \emptyset$. This ensures that the trajectory will always be reset to a point with positive distance from the subsequent guard.
    \item The set of times where a solution to our system reaches the guard (called the set of \textit{impact times}) has a limit point (as happens, for instance, in the case of the bouncing ball with coefficient of restitution less than 1). To exclude this type of situation, we require that the set of impact times be closed and discrete, as in \cite{CoGo19}. 
\end{enumerate}

The above two assumptions will be assumed implicitly throughout the remainder of the paper. 

\end{remark}

\subsection{Interpolation on Hybrid Systems}
SIEs appear in a variety of robotic applications, typically with some control parameters in the dynamics. It is common in practice to develop a \textit{trajectory tracking controller}, which causes the dynamics to converge to a desired trajectory exponentially fast from some set of initial conditions. Once such a controller is designed, one need only select the trajectory to track. For example, one may wish to design this trajectory so that it solves some interpolation problem. In general, however, the energy consumed along the task will depend on the particular choice in trajectory. 

In continuous-time dynamical systems, cubic splines are a ubiquitous choice, because they are piecewise smooth and low energy or even energy-optimal interpolants. In SIEs, however, cubic polynomials may intersect a guard along its trajectory by chance and fail to reach the desired boundary condition. This can be rectified, for instance, by using the minimizers of $J$ with an appropriate potential (as discussed in Section \ref{Sec: Obstacle Avoidance}) used to avoid the guard as needed. For this to be possible, we will assume that obstacle avoidance is feasible whenever necessary.


For some $T \in \R, \ s \in \N$, consider a collection of times $0 = t_1 < t_2 <\cdots < t_s = T$ and knot points $\{\ \xi^{\sigma(n)}_n\}$, where $\sigma: \{1,\ldots,s\} \to \mathcal{V}$ such that for each $n \in \{1,\ldots,s\}$, $\xi^{\sigma(n)}_n= (q^{\sigma(n)}_n, v^{\sigma(n)}_n) \in TQ_{\sigma(n)}$ and $q^{\sigma(n)}_n \notin S_{(\sigma(n),j)}$ for any $j \in \mathcal{N}_{\sigma(n)}$.

\textbf{Interpolation Problem:} Find a piecewise smooth curve $\gamma: [0, T] \to Q$ such that $\gamma(t_n) = q_n^{\sigma(n)}$ and $\dot{\gamma}(t_n) = v_n^{\sigma(n)}$.

\textit{Procedure:} We will construct a sequence of piecewise smooth curves, each of which solve some boundary conditions of the form $\gamma_n(t_n) = \xi_n^{\sigma(n)}$ and $\gamma_n(t_{n+1}) = \xi_{n+1}^{\sigma(n+1)}$. We then need only "glue" the curves together to get the desired result. We consider the two operators as follows: For any $\xi, \eta \in TQ_{i}$, denote by $\mathbb{P}^{\tau}(\xi, \eta; \cdot): [0, \tau] \to Q_i$ a minimizer of $J$ (equation \eqref{J}) on $\Omega_{\xi, \eta}^{\tau}$ corresponding to some potential used to avoid $\cup_{j \in \mathcal{N}_i} S_{ik}$ (as, for example, was designed in Section \ref{sec: large-obstacles}). Denote by $\mathbb{P}^{\tau}_{ij}(\xi, \eta; \cdot): [0, \tau] \to Q_i$ a minimizer corresponding to some potential used to avoid $\cup_{j \in \mathcal{N}_i, k \ne j} S_{ik}$. In other words, $\mathbb{P}^{\tau}$ avoids all guards in $Q_i$, while $\mathbb{P}^{\tau}_{ij}$ avoids all guards except for $S_{ij}$.

\textit{Case 1:} $\sigma(n) = \sigma(n+1)$. In this case, we wish to remain in $Q_{\sigma(n)}$ along the trajectory, and so we must avoid each guard. Clearly it suffices to consider\\ $\gamma_n(t) = \mathbb{P}^{t_{n+1} - t_n}(\xi^{\sigma(n)}_n, \xi^{\sigma(n+1)}_{n+1}; t - t_n)$ for $t \in [t_n, t_{n+1})$.

\textit{Case 2:} $\sigma(n) \ne \sigma(n+1)$. Here, we must move to a different component of the SIEs by passing through the guard(s). In general, it may happen that there is no edge connecting the vertices $y_{\sigma(n)}$ to $y_{\sigma(n+1)}$. However, since the graph is connected by assumption, we can always find a path through the graph connecting $y_{\sigma(n)}$ to $y_{\sigma(n+1)}$. We denote such a path by $y_{\sigma(n)} = y_1^\ast \to y_2^\ast \to...\to y_m^\ast = y_{\sigma(n+1)}$, and further define the edge $S_{i^\ast} := S_{(y_i^\ast, y_{i+1}^\ast)}$ and similarly for $\Delta_{i^\ast}$ and $\mathbb{P}^{\tau}_{i^\ast}$. Set $\alpha_n = \frac{t_{n+1} - t_n}{m}$ and choose a point $\eta_1 \in S_{1^\ast}$. We consider the curve $\gamma_n^1(t) := \mathbb{P}^{\alpha_n}_{1^\ast}(\xi_{\sigma(n)}, \eta_1; t - t_n)$ for $t \in [t_n, t_n + \alpha_n]$. Let $t_n + \tau_1$ denote the first time for which $\gamma(t_n + \tau_1) \in S_{1^\ast}$ (it may happen by chance that the curve intersects the guard before reaching $\eta_1)$. Denote this point of contact by $\xi_1^\ast := \gamma(t_n + \tau_1)$.

We now define the curve $\gamma_n^j$ recursively for $j = 2,...,m-1$. Choose $\eta_j \in S_{j^\ast}$, and let $\tau_j^\ast := t_n + \sum_{k=1}^{j-1} \tau_k$. Define $\gamma_n^j(t) := \mathbb{P}_{j^\ast}^{\alpha_n}(\Delta_{j^\ast}(\xi_{j-1}^\ast), \eta_j; t - \tau_j^\ast)$ for $t \in [\tau_j^\ast, t_n + j\alpha_m]$. Further let $\tau_j$ and $\xi_j^\ast$ be defined by the relation $\xi_j^\ast := \gamma_n^j(\tau_j^\ast + \tau_j) \in S_{j^\ast}$ and $\gamma_n^j(t) \notin S_{j^\ast}$ for any $t < \tau_j^\ast + \tau_j$ in its domain. 

For $j = m,$ we choose \\ $\gamma_n^m(t) := \tilde{\mathbb{P}}^{t_{n+1} - \tau_{m-1}^\ast}(\Delta_{m-1^\ast}(\xi_{m-1}^\ast), \xi^{n+1}_{\sigma(n+1)}; t - \tau_{m-1}^\ast)$\\ for $t \in [\tau_{m-1}^\ast, t_{n+1}]$. Finally, we set
\setlength{\abovedisplayskip}{3pt}
\begin{align*}\gamma_n(t) = \begin{cases}
\gamma_n^1(t) & \text{for } \ t \in [t_n, \tau_1^\ast) \\
\gamma_n^2(t) & \text{for } \ t \in [\tau_1^\ast, \tau_2^\ast)\\
\quad \vdots \\
\gamma_n^m(t) & \text{for } \ t \in [\tau_{m-1}^\ast, t_{n+1}]
\end{cases}
\end{align*}
From which it is clear that $\gamma_n$ is a piecewise smooth curve satisfying $\gamma_n(t_n) = q_n^{\sigma(n)}, \  \dot{\gamma}_n(t_n) = v_n^{\sigma(n)}, \ \gamma_n(t_{n+1}) = q_{n+1}^{\sigma({n+1})},\  \dot{\gamma}_{n+1}(t_{n+1}) = v_{n+1}^{\sigma({n+1})}$ for $n = 1,..., s$.

This strategy, when feasible, produces a solution to the interpolation problem. However, it is worth noting that the procedure (particularly along Case 2) requires some choice in point on the guard. This choice can be made arbitrarily, but one may desire to understand if an ``optimal" choice can be made. We can find necessary conditions for the final point by redoing the derivation of Proposition \ref{th1}, but without the assumption that the variational vector fields vanishing at the end point. The effect of this is that the boundary terms picked up when integrating by parts will no longer vanish at $t = T$, from which necessary conditions can be derived with appropriate choices in the variational vector field. In particular, we have

\begin{lemma}
Suppose $S \subset Q$. A critical point $q$ of $J$ satisfying the initial conditions $q(0) = q_0$ and $\dot{q}(0) = v_0$, and the unilateral constraint $q(T) \in S$ is a smooth curve satisfying
\begin{align}
    \frac{D^3}{dt^3}\dot{q}(t)+R\Big{(}\frac{D}{dt}\dot{q}(t),\dot{q}(t)\Big{)}\dot{q}(t)=- \hbox{\grad} \, V(q(t)),\label{nec} \\
    \Big{|}\Big{|}\frac{D}{dt} \dot{q}(T)\Big{|}\Big{|}_{q(T)}^2 = g_{q(T)}\Big{(}\dot{q}(T), \frac{D^2}{dt^2} \dot{q}(T)\Big{)}.\label{boundary}
\end{align}
\end{lemma}

We see that (\ref{nec}) is identical to the necessary conditions (\ref{eqq1}) where the end point was fixed. On the other hand, (\ref{boundary}) gives us a condition on the derivatives of $q$ at $t = T$. Note, however, that when closed forms for the solutions to (\ref{nec}) cannot be found (which in general cannot be expected), solving (\ref{boundary}) will be very difficult. Moreover, it is not guaranteed that the the trajectory will not intersect the guard at some point before reaching the desired final point (though we may get arbitrarily close by choosing a point near the desired point that is not on the guard and then solving the boundary value problem with the avoidance potential—however this may lead to a very costly trajectory). Note also that condition (\ref{boundary}) is different than the one imposed for velocities to reach a final submanifold in \cite{mishal}.

\begin{remark}
Analogous to  continuous-time dynamical systems, periodic orbits can be defined in a SIEs by a continuous curve $\gamma: \R \to Q$ and a real number $T \in R$ such that $\gamma(t + T) = \gamma(t)$ for all $t \in \R$. The prevailing difference in a SIEs is that such orbits will in general not be closed. Note that the solution to the interpolation problem can be used to construct periodic orbits in a SIEs by choosing the knot points such that $\xi_1^{\sigma(1)} = \xi_s^{\sigma(s)}$.
\end{remark}

\section{Conclusions and Future Work}
In this paper, we presented a motion planning strategy based on variational principles and an artificial potential function. The global existence of extrema was proven in the case that the potential is non-negative. The obstacle avoidance task was studied, and conditions on the artificial potential were derived to guaranteed obstacle avoidance within some tolerance in the case of point-obstacles and obstacles represented by totally bounded subsets of the underlying manifold. Furthermore, a smooth family of potential was provided which may be used to yield obstacle avoidance of an arbitrarily large tolerance (bounded by the geometry of the manifold and the sensing radius). Finally, the results were applied to solve an interpolation problem on Systems with Impulse Effects (SIEs).

Some future work to consider is the application of the results to Lie groups and symmetric spaces, which offers sufficiently rich geometry and symmetry to reduce the necessary conditions \eqref{eqq1}, and may further be used to simplify equation \eqref{boundary}—making the results more directly applicable in applications. We are further interested in studying sufficient conditions for optimality through the second variation of the action \eqref{J} and bi-Jacobi fields.

\section*{Acknowledgements}
Both authors conduct their research at Instituto de Ciencias Matematicas (CSIC-UAM-UC3M-UCM), Calle Nicolas Cabrera 13-15, 28049, Madrid, Spain. The projects that gave rise to these results received the support of a fellowship from ”la Caixa” Foundation (ID 100010434). The fellowship codes are LCF/BQ/DI19/11730028 for Jacob R. Goodman (\textcolor{blue}{jacob.goodman@icmat.es}) and LCF/BQ/PI19/11690016 for Leonardo J. Colombo \\(\textcolor{blue}{leo.colombo@icmat.es}). The authors were also partially funded by Ministerio de Economia, Industria y Competitividad (MINECO, Spain) under grant MTM2016-76702-P  and
``Severo Ochoa Programme for Centres of Excellence'' in
R$\&$D (SEV-2015-0554). All the results are original and have been accepted for publication in the conference proceedings of the 7th IFAC Workshop on Lagrangian and Hamiltonian Methods for Nonlinear Control (LHMNC21) at the Technical University of Berlin.


\begin{thebibliography}{20}
\providecommand{\newblock}{\relax}
\bibitem{mishal}
M. Assif, R. Banavar, A. Bloch, M. Camarinha, L. Colombo.  Variational collision avoidance problems on Riemannian manifolds. in Proceedings of the IEEE International
Conference on Decision and Control, 2018, pp. 2791-2796.

\bibitem{BlCaCoCDC} A. Bloch, M. Camarinha, L. Colombo. Variational obstacle avoidance
on Riemannian manifolds. in Proceedings of the IEEE International
Conference on Decision and Control, 2017, pp. 146-150.
\bibitem{BlCaCoIJC} A. Bloch, M. Camarinha, L. J. Colombo.  Dynamic interpolation for obstacle avoidance on Riemannian manifolds. International Journal of Control, 94 (3), 588-600, 2021. 

\bibitem{point} A. Bloch, M. Camarinha, L. Colombo. Variational point-obstacle avoidance on Riemannian manifolds. Mathematics of Control, Signals, and Systems, 33(1), 109-121, 2021.

\bibitem{Boothby} W.  M.  Boothby. An  Introduction  to  Differentiable  Manifolds  and Riemannian Geometry. Orlando,  FL: Academic  Press Inc., 1975.



\bibitem{marg} M. Camarinha, F. Silva Leite, and P.Crouch. Splines of class $C^k$ on non-euclidean spaces. IMA Journal of
Mathematical Control \& Information, 12:399-410, 1995.
\bibitem{sh} R.S. Chandrasekaran, L. Colombo, M. Camarinha, R. Banavar, A. Bloch. Variational collision and obstacle avoidance of multi-agent systems on Riemannian manifolds. 2020 European Control Conference (ECC), 1689-1694, IEEE, 2020.
\bibitem{CoGo20} L. Colombo, and J. Goodman. A Decentralized Strategy for Variational Collision Avoidance on Complete Riemannian Manifolds. Portuguese Conference on Automatic Control. Springer, Cham, 2020.

\bibitem{CLACC} P. Crouch and F. Silva Leite. Geometry and the Dynamic Interpolation
Problem. Proc. American Control Conference, 1131-1137, 1991.
\bibitem{CroSil:95} P. Crouch, F. Silva Leite. The dynamic interpolation problem: on Riemannian manifolds, Lie groups, and symmetric spaces,   J. Dynam. Control Systems. 1 (1995), no. 2, 177--202.
\bibitem{Giambo} R. Giambò,
F. Giannoni, P. Piccione. An analytical theory for Riemannian cubic polynomials. IMA J. Math Control Information 19:445-460, 2002.
\bibitem{RiemannianPoly} R. Giambò, F. Giannoni, P. Piccione. Optimal Control on Riemannian Manifolds by Interpolation. MCSS 16:278-296, 2004.
\bibitem{teelSurvey}
R. Goebel, R. Sanfelice, and A. Teel. Hybrid dynamical systems.
 Princeton University Press. 2012.
 \bibitem{CoGo19} J. Goodman, L. Colombo. \textit{On the Existence and Uniqueness of Poincar\'e Maps for Systems with Impulse Effects}. IEEE Transactions on Automatic Control $65(4)$, 1815-1821, 2019.
\bibitem{CollAvoid} J. Goodman, and L. Colombo. Variational Collision Avoidance on Riemannian Manifolds. ArXiv preprint arXiv:2104.04285, 2021.
\bibitem{cha} W. Haddad, V. Chellaboina, and S. Nersesov. Impulsive and hybrid dynamical systems. Princeton University Press, 2006.
\bibitem{kod} D. E. Koditschek and E. Rimon. Robot navigation functions on man- ifolds with boundary. Advances in Applied Mathematics, 11(4):412– 442, 1990.
\bibitem{Palais}R.  Palais and Ch.-L. Terng (1988) \textit{Critical Point Theory and Submanifold Geometry}. Berlin: Springer.

\bibitem{hilbert} M. Bauer, C. Maor, and P. W. Michor. Sobolev metrics on spaces of manifold valued curves. arXiv preprint 	arXiv:2007.13315, 2020.

\bibitem{noakes} L. Noakes, G. Heinzinger, B. Paden. Cubic splines on curved spaces. IMA Journal of Mathematical Control and Information, 6(4), 465-473, 1989.
\bibitem{popei} T. Popiel. Higher order geodesics in Lie groups. Math. Control Signals Syst. 19, 235–253 (2007). 

\bibitem{elastica} P. Schrader. Morse theory for elastica. Journal of Geometric Mechanics, 8(2), p.235, 2016.

\end{thebibliography}
\end{document}